\documentclass[a4paper,10pt]{amsart}

\usepackage[english,francais]{babel}
\usepackage{amsmath,enumerate, amsfonts, amssymb,amsthm}
\input xy
\xyoption{all}

\newtheorem{Def}{DEFINITION}[section]
\newtheorem{Theo}[Def]{THEOREME}
\newtheorem{Prop}[Def]{PROPOSITION}
\newtheorem{Rem}[Def]{REMARQUE}
\newtheorem{Lem}[Def]{LEMME}
\newtheorem{Cor}[Def]{COROLLAIRE}
\newtheorem{Ex}[Def]{EXEMPLE}

\newtheorem*{Theo2}{THEOREME}

\DeclareMathOperator{\hess}{hess}

\title[Structures de Monge-Amp\`ere en dimension 6]{Structures de
  Monge-Amp\`ere symplectiques non d\'eg\'en\'er\'ees en dimension 6}

\author{Bertrand BANOS}
\address{Bertrand Banos\\D\'epartement de Math\'ematiques\\
  Universit\'e d'Angers\\ 2 bd Lavoisier, 49045 Angers, France}
\email{bertrand.banos@univ-angers.fr}

\begin{document}

\begin{abstract}
We define a non-degenerated Monge-Amp\`ere structure on a $6$-mani\-fold
associated with a Monge-Amp\`ere equation
as a couple $(\Omega,\omega)$, such that  $\Omega$ is a symplectic form and
$\omega$ is a  $3$-differential form which satisfies $\omega\wedge\Omega=0$
and which is non-degenerated in the sense of Hitchin. We associate with
such a couple an almost (pseudo) Calabi-Yau structure and we study its
integrability from the point of view of Monge-Amp\`ere operators
theory. The result we prove appears as an  analogue of Lychagin and
Roubtsov theorem on integrability of the almost complex or almost product
structure associated with an elliptic or hyperbolic Monge-Amp\`ere
equation in the dimension $4$. We study from this point of view the
example of the Stenzel metric on $T^*S^3$. 

\end{abstract}

\maketitle

\section{Introduction}

Une \'equation de Monge-Amp\`ere est une \'equation diff\'erentielle
du second ordre non lin\'eaire, dont la non-lin\'earit\'e est tr\`es
sp\'ecifique: c'est celle du d\'eterminant. Par exemple en deux
variables, une \'equation de Monge-Amp\`ere s'\'ecrit
\begin{equation}{\label{EMA}}
A\frac{\partial^2 f}{\partial q_1^2}+2B\frac{\partial^2 f}{\partial
  q_1\partial q_2}+ C\frac{\partial^2 f}{\partial q_2^2}+
D(\frac{\partial^2 f}{\partial q_1^2}\frac{\partial^2 f}{\partial
  q_2^2}-(\frac{\partial^2 f}{\partial q_1\partial q_2})^2)+E=0,
\end{equation}
o\`u $A,B,C,D$ et $E$ sont des fonctions lisses sur l'espace des jets
$J^1\mathbb{R}^2$, i.e. sont fonctions de  $(q,f(q),\frac{\partial
  f}{\partial q})$.  Lorsque $A,B,C,D$ et $E$ sont fonctions
uniquement de $(q,\frac{\partial f}{\partial q})$ (i.e. sont des
fonctions sur l'espace cotangent $T^*\mathbb{R}^2$), on parle
d'\'equation de Monge-Amp\`ere symplectique.
 
La th\'eorie g\'eom\'etrique des op\'erateurs de Monge-Amp\`ere
developp\'ee dans les ann\'ees 70 par Lychagin (\cite{L}),
associe aux \'equations de Monge-Amp\`ere symplectiques
sur une vari\'et\'e $M^n$ de dimension $n$ un couple de formes
diff\'erentielles $(\Omega,\omega)$ sur le fibr\'e cotangent $T^*M$,
$\Omega$ \'etant la forme symplectique sur $T^*M$ et $\omega\in
\Omega^n(T^*M)$ \'etant  effective i.e. $\omega\wedge \Omega=0$. Plus
pr\'ecis\'ement, l'\'equation de Monge-Amp\`ere associ\'ee \`a un tel
couple $(\Omega,\omega)$ est l'\'equation diff\'erentielle
\begin{equation}{\label{OMA}}
(df)^*(\omega)=0,
\end{equation}
o\`u $df:M\rightarrow T^*M$ est la section naturelle associ\'ee \`a
une fonction lisse $f$ sur $M$. Ainsi, par exemple, \`a l'\'equation
$\eqref{EMA}$ est associ\'ee la forme symplectique canonique sur
$T^*\mathbb{R}^2$ 
$$
\Omega= dq_1\wedge dp_1+dq_2\wedge dp_2,
$$
et la forme diff\'erentielle
$$
\begin{aligned}
\omega&=A dp_1\wedge dq_2 + B(dq_1\wedge dp_1-dq_2\wedge dp_2)+C
dq_1\wedge dp_2\\
&+D dp_1\wedge dp_2+ E dq_1\wedge dq_2.\\
\end{aligned}
$$
Une solution g\'en\'eralis\'ee de l'\'equation de Monge-Amp\`ere
$\eqref{OMA}$  est une sous vari\'et\'e lagrangienne de
$(T^*M,\Omega)$ sur laquelle s'annule notre forme $\omega$. Remarquons
que si une solution g\'en\'eralis\'ee se projette diff\'eomorphiquement
sur la base $M$ alors cette solution est localement le graphe d'une
section $df: M\rightarrow T^*M$ avec $f$ solution explicite de
$\eqref{OMA}$.

G\'en\'eralisant cette notion, nous appelons structure de
Monge-Amp\`ere (symplectique) sur une vari\'et\'e $X^{2n}$ de dimension $2n$ la
donn\'ee d'un couple $(\Omega,\omega)$, $\Omega$ \'etant une forme
symplectique sur $X$ et $\omega\in \Omega^n(X)$ \'etant une forme
diff\'erentielle effective sur $(X,\Omega)$. En dimension $n=2$,
lorsque la $2$-forme $\omega$ est non d\'eg\'en\'er\'ee (i.e. le
pfaffien $pf(\omega)=\frac{\omega\wedge\omega}{\Omega\wedge\Omega}$
est non nul), le couple $(\Omega,\omega)$ d\'efinit une section
$A_\omega:X\rightarrow TX\otimes T^*X$  qui est soit une structure presque complexe
($A_\omega^2=-Id$) lorsque $pf(\omega)>0$, soit une
structure presque produit $(A_\omega^2=Id)$ lorsque $pf(\omega)<0$:
$$
\frac{\omega}{\sqrt{|pf(\omega)|}}=\Omega(A_\omega.,.).
$$
Lychagin et Roubtsov (\cite{LR1}) ont montr\'e la 
\begin{Prop}{\label{int2}}
La structure presque complexe ou presque produit $A_\omega$ est
int\'egrable si et seulement si  
$$
d(\frac{\omega}{\sqrt{|pf(\omega)|}})=0.
$$
\end{Prop}

Dans le langage des \'equations diff\'erentielles, c'est \`a dire du
point de vue local, ce r\'esultat s'\'enonce
\begin{Prop}{\label{class2}}
Une \'equation de Monge-Amp\`ere $\Delta_\omega=0$ sur $\mathbb{R}^2$
associ\'ee \`a une forme effective non d\'eg\'en\'er\'ee $\omega$ 
peut \^etre ramen\'ee  par un changement de variable symplectique \`a l'une
des deux \'equations
$$
\begin{cases}
\Delta f=0&(pf(\omega)>0)\\
\square f=0&(pf(\omega)<0)\\
\end{cases}
$$
si et seulement si 
$$
d(\frac{\omega}{\sqrt{|pf(\omega)|}})=0.
$$
\end{Prop}

Pour classifier les \'equations de Monge-Amp\`ere \`a coefficients
constants en dimension $3$, Lychagin et Roubtsov ont introduit un invariant
quadratique $q_\omega$ associ\'e \`a chaque $3$-forme effective $\omega$
(\cite{LR3}). Hitchin, en \'etudiant les structures de Calabi-Yau
en dimension complexe $3$, a d\'efini un invariant lin\'eaire
$K_\omega$ (\cite{Hi}) pour toutes les $3$-formes. Nous montrons ici
que ces deux invariants co\"\i ncident lorsque l'on se restreint aux formes
effectives. (Proposition \ref{compatibilite}). Cette remarque nous permet de
concilier les travaux  de Lychagin, Roubtsov et Hitchin  pour
d\'emontrer des r\'esultats  analogues \`a $\ref{int2}$ et
$\ref{class2}$  en dimension $3$. Pr\'ecisons ici
quelques notations.  Nous supposons que notre $3$-forme effective
$\omega$ est non d\'eg\'en\'er\'ee au sens de Hitchin. Hitchin a
d\'efini un invariant scalaire $\lambda(\omega)$ que nous appelerons pfaffien
de Hitchin et qui est alors non nul. Il a de plus montr\'e que
$\omega$  s'\'ecrit comme la somme de deux formes
complexes d\'ecomposables, uniquement d\'etermin\'ees \`a ordre
pr\`es: $\omega=\alpha+\beta$, et que cette d\'ecomposition permet
d'associer une forme duale $\hat{\omega}$ \`a $\omega$. Nous
associons au  couple $(\Omega,\omega)\in \Omega^2(X^6)\times
\Omega^3(X^6)$   une structure g\'eom\'etrique que nous appelons
structure presque (pseudo) Calabi-Yau. Cette structure est
essentiellement la donn\'ee d'une m\'etrique $q_\omega$
(\'eventuellement ind\'efinie), d'une structure presque complexe ou
presque produit $K_\omega$ compatible avec $q_\omega$ et $\Omega$ et
de deux formes diff\'erentielles de degr\'e $3$ d\'ecomposables dont
les distributions associ\'ees sont les distributions des  vecteurs
propres de $K_\omega$.  Nous \'etendons les r\'esultats de \cite{Hi}
pour montrer l'analogue de $\ref{int2}$:

\begin{Prop}
La structure presque (pseudo) Calabi-Yau
$(q_\omega,K_\omega,\Omega,\alpha,\beta)$ est ``int\'egrable'' si et
seulement si   
$$
\begin{cases}
d(\frac{\omega}{\sqrt[4]{|\lambda(\omega)|}})=0&\\
d(\frac{\hat{\omega}}{\sqrt[4]{|\lambda(\omega)|}})=0&\\
\end{cases}
$$
\end{Prop}

Nous montrons alors la version locale de cette proposition:

\begin{Theo2}{\bf\ref{theo}}
Une \'equation de Monge-Amp\`ere symplectique en dimension $3$
associ\'ee \`a une structure de Monge-Amp\`ere non d\'eg\'en\'er\'ee
$(\Omega,\omega)$ peut \^etre ramen\'ee par un changement de variables
symplectique \`a l'une des \'equations suivantes
$$
\begin{cases}
\hess (f)=\gamma,&\gamma\neq 0\\
\Delta f-\gamma \hess (f)=0, &\gamma\neq 0\\
\square f+\gamma \hess (f)=0, &\gamma \neq 0\\
\end{cases}
$$
si et seulement si
$$
\begin{cases}
d(\frac{\omega}{\sqrt[4]{|\lambda(\omega)|}})=0&\\
d(\frac{\hat{\omega}}{\sqrt[4]{|\lambda(\omega)|}})=0&\\
q_\omega\text{ est plate}&\\
\end{cases}
$$
\end{Theo2}

La motivation de cette article est de g\'en\'eraliser en dimension $3$ la
notion de ``sous vari\'et\'e sp\'eciale lagrangienne''. Ces
sous-vari\'et\'es  ont \'et\'e introduites dans les ann\'ees 80 par
Harvey et Lawson dans leur c\'el\`ebre article
\cite{HL}. Initialement, ils recherchaient des exemples de sous
vari\'et\'es minimales. Rappelons qu'une $p$-calibration sur une
vari\'et\'e riemannienne $(Y,g)$ est une $p$-forme ferm\'ee $\phi\in
\Omega^p(Y)$ telle que en tout point $y$ de $Y$ 
$$
|\phi_y(e_1,\ldots,e_p)|\leq 1,
$$
pour toute famille orthonorm\'ee $(e_1,\ldots,e_p)$ de $T_yY$. Une
sous-vari\'et\'e orient\'ee $L$ de $Y$ de dimension $p$ est dite
$\phi$-calibr\'ee si $\phi_y(\theta_{L,y})=1$ pour tout $y\in
L$, o\`u $\theta_{L,y}$ est la forme volume sur $T_yL$ d\'efinie par la
m\'etrique $g$ et l'orientation choisie sur $L$.  Les
sous-vari\'et\'es $\phi$-calibr\'ees sont alors de volume minimal dans
leur classe d'homologie. La forme
$Re(\alpha)$ est un exemple de calibration sur
$\mathbb{C}^n$ avec $\alpha=dz_1\wedge \ldots \wedge dz_n$ et les
sous-vari\'et\'es sp\'eciales lagrangiennes sont les sous vari\'et\'es
$Re(\alpha)$-calibr\'ees. Cette notion de calibration sp\'eciale
lagrangienne se g\'en\'eralise sur les vari\'et\'es de Calabi-Yau,
i.e. sur les vari\'et\'es de K\"ahler $(Y,g,I,\Omega)$ de dimension
complexe $n$ munies d'une forme volume holomorphe $\alpha$ telle que
$\frac{\alpha\wedge\overline{\alpha}}{\Omega^n}$ soit une fonction
constante non nulle.   
Comme l'explique M. Audin dans son cours \emph{Lagrangian submanifold}
(\cite{A}), nous avons pu assister
ces derni\`eres ann\'ees \`a un net regain d'int\'er\^et pour ces sous
vari\'et\'es calibr\'ees apr\`es les travaux de Strominger, Yau et
Zaslow en ``th\'eorie miroir'' (\cite{SYZ}) et notamment leur
construction d'une vari\'et\'e ``miroir'' d'une vari\'et\'e de
Calabi-Yau de dimension complexe $3$ \`a partir d'une fibration
torique sp\'eciale lagrangienne.  

En 1990, Gromov (dans une discussion avec Roubtsov) avait not\'e que
les structures de Monge-Amp\`ere sont, en un certain sens,
l'\'equivalent symplectique des calibrations: les calibrations
correspondent aux formes effectives et les sous vari\'et\'es
calibr\'ees correspondent aux sous vari\'et\'es lagrangiennes. Les
sous-vari\'et\'es sp\'eciales lagrangiennes se trouvent \`a
l'intersection des deux approches.  Harvey et Lawson ont en effet
montr\'e, que les sous-vari\'et\'es sp\'eciales lagrangiennes de
$\mathbb{C}^n$ sont, \`a choix d'orientation pr\`es, les solutions
g\'en\'eralis\'ees  de l'\'equation diff\'erentielle associ\'ee \`a la
structure de  Monge-Amp\`ere $(\Omega,\omega)$ avec 
$$
\begin{cases}
\Omega=\frac{i}{2}(dz_1\wedge d\overline{z_1}+\ldots + dz_n\wedge
d\overline{z_n})&\\
\omega=Im(dz_1\wedge \ldots \wedge dz_n)&\\
\end{cases}
$$
Cette \'equation s'\'ecrit en dimension $3$  
\begin{equation}{\label{slag}} 
\Delta f-\hess (f)=0.
\end{equation}

Nous voulons montrer ici comment l'on peut  construire des structures
g\'eom\'etriques, analogues du point de vue symplectique  aux
structures de Calabi-Yau classiques, \`a partir d'\'equations
diff\'erentielles semblables \`a $\eqref{slag}$.

Dans une premi\`ere partie, nous rappelons l'approche g\'eom\'etrique de
Lychagin des \'equations de Monge-Amp\`ere \`a partir des formes
diff\'erentielles effectives. Nous illustrons cette approche sur un
exemple issu du mod\`ele ``semi-g\'eostrophique'': les \'equations de
Chynoweth-Sewell. Nous remarquons que ces \'equations sont
symplectiquement \'equivalentes \`a l'\'equation de Monge-Amp\`ere
classique $\hess (f)=1$ en trois variables  et nous donnons une solution
explicite relativement ``g\'en\'erique'' de cette \'equation.  Dans une deuxi\`eme
partie nous adaptons les travaux de Hitchin sur la g\'eom\'etrie des
$3$-formes diff\'erentielles pour les formes effectives. Enfin, dans
une troisi\`eme partie, nous d\'efinissons les structures presque
Calabi-Yau associ\'ees aux structures de Monge-Amp\`ere, nous
\'etudions leur int\'egrabilit\'e et nous montrons le th\'eor\`eme
\ref{theo}. Nous \'etudions comme exemple de m\'etrique de Calabi-Yau
non plate,  la m\'etrique de Stenzel sur $T^*S^3$.

Cet article constitue une partie de la th\`ese
 de l'auteur pr\'epar\'ee \`a l'universit\'e d'Angers. Je voudrais
 remercier mon directeur de th\`ese V. Roubtsov pour m'avoir
 sugg\'erer et expliquer ce probl\`eme et pour les tr\`es
 enrichissantes   discussions que nous avons pu avoir. Je remercie mon
 coll\`egue Oleg Lisovyy pour toute l'aide qu'il m'a apport\'ee,
 notamment dans la construction d'une solution explicite
 des \'equations de Chynoweth-Sewell. Je voudrais aussi  remercier le
 Professeur B. Enriquez pour toutes ses remarques et suggestions.

\section{Formes effectives et structures de Monge-Amp\`ere}

Soit $(V,\Omega)$ un espace vectoriel symplectique r\'eel de dimension $2n$. Soit
$\Gamma:V\rightarrow V^*$ l'isomorphisme induit par $\Omega$ et $X_\Omega\in
\Lambda^2(V^*)$ l'unique bivecteur tel que $\Gamma^*(X_\Omega)=\Omega$,
$\Gamma^*:\Lambda^*(V)\rightarrow \Lambda^*(V^*)$ \'etant la puissance ext\'erieure de
$\Gamma$ \footnote{on note par $\Lambda^*(V^*)$ l'espace des formes ext\'erieures sur
un espace vectoriel $V$ et $\Omega^*(X)$ l'espace des formes diff\'erentielles sur une
vari\'et\'e lisse $X$}.

Suivant les notations de (\cite{L}), nous introduisons les op\'erateurs $\bot:
\Lambda^k(V^*)\rightarrow \Lambda^{k-2}(V^*)$ $\omega\mapsto i_{X_\Omega}(\omega)$ et
$\top: \Lambda^k(V^*)\rightarrow \Lambda^{k+2}(V^*)$, $\omega\mapsto \omega\wedge
\Omega$. Ces op\'erateurs v\'erifient:
$$
\begin{cases}
[\bot,\top](\omega)=(n-k)\omega,\;\forall\omega\in\Lambda^k(V^*);&\\
\bot:\Lambda^k(V^*)\rightarrow \Lambda^{k-2}(V^*)\text{ est injective pour $k\geq n+1$};&\\
\top:\Lambda^k(V^*)\rightarrow \Lambda^{k+2}(V^*)\text{ est injective pour $k\leq n-1$}.&\\
\end{cases}
$$

Une $k$-forme $\omega$ est dite effective si $\bot\omega=0$. Nous noterons
$\Lambda^k_\varepsilon(V^*)$ l'espace des $k$-formes effectives. Si $k=n$, $\omega$ est
effective si et seulement si $\omega\wedge \Omega=0$.

Le th\'eor\`eme suivant explique le r\^ole fondamental jou\'e par les formes effectives
dans la th\'eorie des op\'erateurs de Monge-Amp\`ere (\cite{L}):

\begin{Theo}[Hodge-Lepage-Lychagin]{\label{Hodge}}
\begin{enumerate}
\item Toute forme $\omega\in \Lambda^k(V^*)$ peut \^etre d\'ecompos\'ee de fa\c con
unique en une somme finie
$$
\Omega= \omega_0+\top\omega_1+\top^2\omega_2+\ldots,
$$
les formes $\omega_i$ \'etant toutes effectives.
\item Si deux $k$-formes effectives s'annulent sur les m\^eme sous espaces isotropes de
dimensions $k$ de $(V,\Omega)$ alors elles sont proportionnelles.
\end{enumerate}
\end{Theo}

Soit $M^n$ une vari\'et\'e lisse de dimension $n$. Notons $J^1M$ le fibr\'e des
$1$-jets des fonctions lisses sur $M$ et notons $j^1(f): M\rightarrow J^1M$, $x\mapsto
[f]^1_x$ la section naturelle associ\'ee \`a une fonction lisse
$f$. L'op\'erateur de  Monge-Amp\`ere $\Delta_\omega:
C^\infty(M)\rightarrow \Omega^n(M)$ associ\'e \`a une forme $\omega\in
\Omega^n(J^1M)$ est  d\'efini par
$$
\Delta_\omega(f)= j^1(f)^*(\omega).
$$

Soit $U$ la forme de contact sur $J^1M$ et $X_1$ le champ de Reeb. Une solution
g\'en\'eralis\'ee de l'\'equation de Monge-Amp\`ere $\Delta_\omega=0$ est une sous
vari\'et\'e legendrienne $L^n$ de $(J^1M,U)$ sur laquelle s'annule $\omega$. Si $L$ se
projette diff\'eomorphiquement sur $M$ alors $L$ est localement le graphe d'une section
$j^1(f)$, $f$ \'etant une solution de l'\'equation diff\'erentielle
$\Delta_\omega(f)=0$. Soit la distribution de contact $x\in
J^1M\mapsto C(x)= Ker(U_x)$. Remarquons
que  $(C(x),dU_x)$ est un espace vectoriel symplectique de dimension
$2n$ et que
$$
T_xJ^1M=C(x)\oplus X_{1x}.
$$
De plus, si $L$ est une sous vari\'et\'e legendrienne alors $T_xL$ est un sous espace
lagrangien de $(C(x),dU_x)$.

Soit $\Omega^*(C^*)$ l'espace des formes s'annulant le long de $X_1$. En tout point
$x$, l'espace vectoriel $(\Omega^*(C^*))_x$ s'identifie naturellement
avec $\Omega^*(C(x)^*)$. Nous noterons $\Omega^*_\varepsilon(C^*)$
l'espace des formes effectives sur $(C(x),dU_x)$ en tout point $x$ de
$J^1M$.  D'apr\`es la premi\`ere partie du th\'eor\`eme \ref{Hodge}, on a
$$
\Omega^*_\varepsilon(C^*)=\Omega^*(J^1M)/I_C,
$$
$I_C$ \'etant l'id\'eal de Cartan engendr\'e par $U$ et $dU$. La seconde partie dit que
si deux formes diff\'erentielles $\omega$ et $\theta$ sur $J^1M$ d\'eterminent le
m\^eme op\'erateur de Monge-Amp\`ere alors $\omega-\theta\in I_C$.

Le pseudo-groupe des diff\'eomorphismes de contacts sur $J^1M$ agit
sur ces
op\'erateurs comme suit
$$
F(\Delta_\omega)= \Delta_{F^*\omega},
$$
et l'action infinit\'esimale associ\'ee est
$$
X(\Delta_\omega)=\Delta_{L_X(\omega)}.
$$

Nous nous int\'eresserons ici \`a une classe plus restrictive d'op\'erateurs, celles des
op\'erateurs symplectiques, i.e. la classe des op\'erateurs qui v\'erifient
$$
X_1(\Delta_\omega)=\Delta_{L_{X_1}(\omega)}=0.
$$
Ces op\'erateurs correspondent aux \'equations de Monge-Amp\`ere
$\Sigma=\Sigma(x,\frac{\partial f}{\partial
  x},\frac{\partial^2 f}{\partial x^2})\subset J^2M$ ne d\'ependant que
des d\'eriv\'ee partielles du premier et du deuxi\`eme ordre de $f$.
Soit $(T^*M,\Omega)$ le fibr\'e cotangent de $M$ muni de sa forme
symplectique naturelle et consid\'erons la projection $\beta:
J^1M\rightarrow T^*M$ d\'efinie par le diagramme commutatif suivant:
$$
\xymatrix{
  \mathbb{R}&&J^1M\ar[ll]_\alpha\ar[rr]^\beta&&T^*M\\
&&&&\\
&&M\ar[lluu]^f\ar[uu]_{j^1(f)}\ar[rruu]_{df}&&\\
}
$$

En utilisant $\beta$, on peut naturellement identifier l'espace $\{\omega\in
\Omega_\varepsilon(C^*): L_{X_1}\omega=0\}$ avec l'espace des formes effectives sur
$(T^*M,\Omega)$. L'action du groupe de contact sur les op\'erateurs symplectiques est
alors  restreinte \`a l'action du groupe symplectique sur $T^*M$.

\begin{Def}
Une structure de Monge-Amp\`ere symplectique sur une vari\'et\'e $X^{2n}$ de dimension
$2n$ est la donn\'ee d'un couple $(\Omega,\omega)$, o\`u $\Omega\in \Omega^2(X)$ est une
forme symplectique sur $X$ et $\omega\in \Omega^n(X)$ est une forme diff\'erentielle
v\'erifiant
$$
\Omega\wedge \omega=0
$$
\end{Def}

Remarquons que, d'apr\`es le th\'eor\`eme de Darboux, $(X^{2n},\Omega)$ s'identifie
localement avec $(T^*\mathbb{R}^n,\Omega_0)$. Au couple $(\Omega,\omega)$ est alors
associ\'ee l'\'equation de Monge-Amp\`ere (symplectique) sur $\mathbb{R}^n$,
$\Delta_\omega=0$. Et r\'eciproquement \`a toute \'equation de Monge-Amp\`ere
(symplectique) $\Delta_\omega=0$ sur une vari\'et\'e $M$ est associ\'ee une classe
conforme de  structures de Monge-Amp\`ere (symplectiques) $(\Omega,\omega)$ sur $T^*M$.

\begin{Ex}
Un exemple  d'\'equation de Monge-Amp\`ere \`a coefficients constant est celui des
\'equations de Chynoweth-Sewell (\cite{CS}) issues du ``mod\`ele semi-g\'eostro\-phique'':
\begin{equation}{\label{chy-sew}}
\frac{\partial^2 f}{\partial x^2}\frac{\partial^2 f}{\partial
  y^2}-(\frac{\partial^2 f}{\partial x\partial y})^2+\frac{\partial^2
  f}{\partial z^2}=\gamma, \;\; \gamma\in \mathbb{R}
\end{equation}

L'\'equation initiale de Chynoweth-Sewell correspond au cas $\gamma=0$. Remarquons que
si $F(x,y)$ est une solution de $\hess (F)=1$ alors $F(x,y)-\frac{1}{2}z^2$ est une
solution de l'\'equation de Chynoweth-Sewell. Ainsi, par exemple
$$
\frac{1}{3}\sqrt{(x^2+2y)^3} -\frac{1}{2} z^2
$$
est une de solution r\'eguli\`ere de $\eqref{chy-sew}$ quand $\gamma=0$.

La forme effective associ\'ee \`a $\eqref{chy-sew}$  s'\'ecrit dans le syst\`eme de
coordonn\'ees symplectique $(x,y,z,p,q,h)$ de $T^*\mathbb{R}^3$,
$$
\omega= dp\wedge dq\wedge dz+dx\wedge dy\wedge dh-\gamma dx\wedge dy\wedge dz.
$$
$\omega$ est clairement la somme de deux formes d\'ecomposables:
$$
\omega=dp\wedge dq\wedge dz+dx\wedge dy\wedge(dh-\gamma dz)
$$
Soit alors le changement de variable symplectique
$$
\phi(x,y,z,p,q,h)=(x,y,h,p,q,\gamma h -z).
$$
L'image par $\phi$ de $\omega$ est
$$
\phi^*(\omega)= dp\wedge dq\wedge dh - dx\wedge dy\wedge dz.
$$
Autrement dit, les \'equations de Chynoweth-Sewell se rapportent toutes \`a l'\'equation
de Monge-Amp\`ere classique
\begin{equation}{\label{hess}}
\hess (f)=1
\end{equation}

Lorsque $n$ est impair, un exemple de solution r\'eguli\`ere de
$\hess (f)=1$ est
$$
f(x_1,\ldots,x_n)=2\int_{a}^{\sqrt{\underset{i<j}{\sum} x_ix_j}}
(b+\frac{\xi^n}{n-1})^{\frac{1}{n}}d\xi,
$$
o\`u $a$ et $b$ sont des constantes. En particulier pour $n=3$,
$$
f(x,y,z)=\int_{a}^{\sqrt{xy+yz+zx}} (b+4\xi^3)^{1/3}d\xi,
$$
est une solution r\'eguli\`ere de $\eqref{hess}$. D\`es lors un
exemple de solution g\'en\'eralis\'ee de $\eqref{chy-sew}$ est
$$
L=\{(x,y,(x+y)\alpha,(y+z)\alpha,(z+x)\alpha,\gamma(x+y)\alpha-z)\},
$$
avec 
$$
\alpha=
\frac{1}{2}(\frac{b}{(xy+yz+zx)^{\frac{3}{2}}}+4)^{\frac{1}{3}}
$$
\end{Ex}

\section{G\'eom\'etrie des $3$-formes effectives en dimension $6$}

\subsection{Action de $SL(6,\mathbb{R})$ sur $\Lambda^3(\mathbb{R}^6)$.\\}

Nous rappelons ici bri\`evement les r\'esultats de Hitchin sur la
g\'eom\'etrie des $3$-formes ext\'erieures en dimension $6$
(\cite{Hi}). Soit $V$ un espace vectoriel r\'eel de
dimension $6$ et $\Lambda^k(V^*)$ l'espace des $k$-formes
ext\'erieures sur $V$. Fixons $\theta\in \Lambda^6(V^*)$ une forme
volume sur $V$.  Notons $A:\Lambda^5(V^*)\rightarrow V\otimes
\Lambda^6(V^*)$ l'isomorphisme induit par le produit ext\'erieur
et soit pour $\omega\in \Lambda^3(V^*)$ l'endomorphisme $K_\omega:
V\rightarrow V$ d\'efini par
$$
K_\omega(X)\theta= A(i_X\omega\wedge \omega).
$$

\begin{Def}
Le pfaffien de Hitchin d'une $3$-forme $\omega\in \Lambda^3(V^*)$
est
$$
\lambda_\theta(\omega)= \frac{1}{6}Tr(K_\omega^2).
$$
Une $3$-forme $\omega$ est dite non d\'eg\'en\'er\'ee si et seulement si
$\lambda_\theta(\omega)$ est non nul.
\end{Def}

\begin{Prop}[Hitchin]{\label{decompo}}
Soit $\omega\in \Lambda^3(V^*)$ de pfaffien $\lambda_\theta(\omega)$ non nul.
\begin{enumerate}
\item $K_\omega^2=\lambda_\theta(\omega)Id$.
\item $\lambda_\theta(\omega)>0$ si et seulement si $\omega=\alpha+\beta$
o\`u $\alpha$ et $\beta$ sont des formes r\'eelles d\'ecomposables
sur $V$. De plus si on impose  $\frac{\alpha\wedge \beta}{\theta}>0$
alors  $\alpha$ et $\beta$ sont uniquement d\'etermin\'ees:
$$
\begin{cases}
2\alpha=\omega+ |\lambda_\theta(\omega)|^{-\frac{3}{2}}K_\omega^*(\omega)&\\
2\beta=\omega- |\lambda_\theta(\omega)|^{-\frac{3}{2}}K_\omega^*(\omega)&\\
\end{cases}
$$
\item $\lambda_\theta(\omega)<0$ si et seulement si $\omega= \alpha+
\overline{\alpha}$ o\`u $\alpha\in \Lambda^3(V^*\otimes
\mathbb{C})$ est une forme complexe d\'ecomposable sur $V$.  Si l'on
impose de plus que  $\frac{\alpha\wedge \overline{\alpha}}{i\theta}>
0$ alors $\alpha$ est uniquement d\'etermin\'ee:
$$
\alpha= \omega+i|\lambda_\theta(\omega)|^{-\frac{3}{2}}K_\omega^*(\omega).
$$
\end{enumerate}
\end{Prop}

\begin{Rem}
Fixons une base $(e_1,\ldots,e_6)$ de $V$ et notons
$(e_1^*,\ldots,e_6^*)$ sa base duale.
\begin{enumerate}
\item $\lambda_\theta(\omega)>0$ si et seulement si $\omega$ est dans la
$GL(V)$-orbite de
$$
e_1^*\wedge e_2^*\wedge e_3^*+ e_4^*\wedge e_5^*\wedge e_6^*.
$$
\item $\lambda_\theta(\omega)<0$ si et seulement si $\omega$ est dans la
$GL(V)$-orbite de
$$
(e_1^*+ie_4^*)\wedge (e_2^*+ie_5^*)\wedge
(e_3^*+ie_6^*)+(e_1^*-ie_4^*)\wedge (e_2^*-ie_5^*)\wedge
(e_3^*-ie_6^*).
$$
\end{enumerate}
L'action de $GL(V)$ sur $\Lambda^3(V^*)$ \`a ainsi deux orbites ouvertes s\'epar\'ees
par l'hyper\-surface $\lambda_\theta=0$. Ceci justifie la terminologie  de "$3$-forme
non d\'eg\'en\'er\'ee en dimension $6$".
\end{Rem}

L'unicit\'e de la d\'ecomposition en somme de deux formes
d\'ecomposables (\`a choix d'orientation pr\`es) permet d'associer
une forme duale $\hat{\omega}$ \`a toute forme non
d\'eg\'en\'er\'ee $\omega$:
\begin{Def}[Hitchin]
\begin{enumerate}
\item Si $\lambda_\theta(\omega)>0$ et  $\omega=\alpha+\beta$ alors
$\hat{\omega}=\alpha-\beta$.
\item Si $\lambda_\theta(\omega)<0$ et $\omega= \alpha+\overline{\alpha}$
alors $\hat{\omega}= i(\overline{\alpha}-\alpha)$.
\end{enumerate}
\end{Def}

Notons enfin que le produit ext\'erieur d\'efinit une structure
symplectique $\Theta$ sur $\Lambda^3(V^*)$
$$
\Theta(\omega,\omega')\theta=\omega\wedge \omega',
$$
pour laquelle l'action de $SL(6,\mathbb{R})$ est hamiltonienne.

\begin{Prop}[Hitchin]
L'action de $SL(6,\mathbb{R})$ sur $(\Lambda^3(V^*),\Theta)$ est
hamiltonienne d'application moment $K: \Lambda^3(V^*)\rightarrow
sl(6,\mathbb{R})$.
\end{Prop}

\subsection{Action de $Sp(6)$ sur $\Lambda_\varepsilon(\mathbb{R}^6)$.\\}

Supposons maintenant que $V$ est muni d'une forme symplectique $\Omega$. Nous fixons la
forme volume $\theta= -\frac{\Omega^3}{6}$ et nous noterons
$\lambda(\omega)=\lambda_\theta(\omega)$ le pfaffien associ\'e \`a une $3$-forme
$\omega$. Nous dirons qu'une base $(e_1,e_2,e_3,f_1,f_2,f_3)$ est symplectique si
$$
\Omega= e_1^*\wedge f_1^*+ e_2^*\wedge f_2^*+ e_3^*\wedge f_3^*.
$$
Une $3$-forme $\omega$ sur $V$ est effective si
$$
\omega\wedge \Omega=0.
$$
Nous noterons  $\Lambda^3_\varepsilon(V^*)$ l'espace des $3$-formes effectives sur
$V$.  $(\Lambda^3_\varepsilon(V^*),\Theta)$ est un sous espace symplectique (de
dimension $14$) de $(\Lambda^3(V^*),\Theta)$ puis\-que, d'apr\`es
\ref{Hodge} toute $3$-forme $\omega$ sur $V$ peut s'\'ecrire 
$$
\omega=\omega_0+\omega_1\wedge \Omega,
$$
avec $\omega_0$ effective.

\begin{Lem}
$\omega\in \Lambda^3(V^*)$ est effective si et seulement si
$K_\omega\in sp(6)$ o\`u $sp(6)$ est l'alg\`ebre de Lie du groupe
symplectique  $Sp(6)=Sp(\Omega)$.
\end{Lem}

\begin{proof}
Dans une base symplectique $(e_1,e_2,e_3,f_1,f_2,f_3)$, $K_\omega$ s'\'ecrit
$$
K_\omega(X)\theta= \sum_{j=1}^3 (i_X\omega\wedge \omega\wedge
e_j^*)\otimes e_j+ \sum_{j=1}^3 (i_X\omega\wedge \omega\wedge
f_j^*)\otimes f_j.
$$
Ainsi, si on note $K_\omega= \begin{pmatrix}
A&B\\C&D\\\end{pmatrix}$, on a
$$
\begin{cases}
A_{jk}\theta= i_{e_j}\omega\wedge \omega\wedge e_k^*&\\
B_{jk}\theta= i_{f_j}\omega\wedge \omega\wedge e_k^*&\\
C_{jk}\theta= i_{e_j}\omega\wedge \omega\wedge f_k^*&\\
D_{jk}\theta= i_{f_j}\omega\wedge \omega\wedge f_k^*&\\
\end{cases}
$$
Or $\omega$ est effective si et seulement si les relations suivantes
sont v\'erifi\'ees pour $k=1,2,3$,
$$
\begin{cases}
i_{e_k}\omega\wedge \Omega=\omega\wedge
i_{e_k}\Omega=\omega\wedge f_k^*&\\
i_{f_k}\omega\wedge \Omega=\omega\wedge
i_{f_k}\Omega=-\omega\wedge e_k^*&\\
\end{cases}
$$
et donc $\omega$ est effective si et seulement si $D=-A^t$,
$B^t=B$ et $C^t=C$ i.e. si et seulement si $K_\omega\in sp(6)$.
\end{proof}

\begin{Cor}
L'action de $Sp(6)$ sur $(\Lambda^3_\varepsilon(V^*),\Theta)$ est
ha\-mil\-tonien\-ne d'appli\-cation moment $K:
\Lambda^3_\varepsilon(V^*)\rightarrow sp(6)$.
\end{Cor}

Dans \cite{LR2} Lychagin et Roubtsov ont associ\'e \`a chaque forme
effective $\omega\in \Lambda^3_\varepsilon(V^*)$  un invariant
symplectique $q_\omega\in S^2(V^*)$ d\'efini par
$$
q_\omega(X)= -\frac{1}{4} \bot^2(i_X\omega\wedge i_X\omega).
$$
Cette invariant mesure en fait les racines du polyn\^ome
caract\'eristique de $\omega$:
$$
(i_X\omega -\xi\Omega)^3=
-\xi(\xi-\sqrt{|q_\omega(X)|})(\xi+\sqrt{|q_\omega(X)|})\Omega^3.
$$
A partir de cet invariant ils ont \'etabli une classification des
diff\'erentes orbites de l'action de $Sp(6)$ sur
$\Lambda^3_\varepsilon(V^*)$. Cette classification a \'et\'e affin\'ee dans
\cite{Ba} et est r\'esum\'ee dans le tableau \ref{tableau}.
\begin{table}[hbp!]
\begin{tabular}{|c|c|c|c|}

\hline

& $\omega$& $q_\omega$&$\lambda(\omega)$\\

\hline

1.& $e_{1}^*\wedge e_{2}^*\wedge e_{3}^*+ \gamma f_{1}^*\wedge
f_{2}^*\wedge f_{3}^*,\gamma \neq
0$&$\frac{\gamma}{2}(e_{1}^*f_{1}^*+
e_{2}^*f_{2}^*+ e_{3}^*f_{3}^*)$&$\gamma^4$\\

\hline

2.&$f_{1}^*\wedge e_{2}^*\wedge e_{3}^*+ f_{2}^*\wedge
e_{1}^*\wedge e_{3}^*$&
$(e_{1}^*)^2-(e_{2}^*)^2+(e_{3}^*)^2$&$-4\gamma^4$\\
&$+ f_{3}^*\wedge e_{1}^*\wedge e_{2}^*+ \gamma^2 f_{1}^*\wedge
f_{2}^*\wedge f_{3}^*,\gamma\neq
0$&$+\gamma^2((f_{1}^*)^2-(f_{2}^*)^2+(f_{3}^*)^2)$&\\

\hline

3.&$f_{1}^*\wedge e_{2}^*\wedge e_{3}^*- f_{2}^*\wedge
e_{1}^*\wedge e_{3}^*$&
$-(e_{1}^*)^2-(e_{2}^*)^2-(e_{3}^*)^2$&$-4\gamma^4$\\
&$+ f_{3}^*\wedge e_{1}^*\wedge e_{2}^*- \gamma^2 f_{1}^*\wedge
f_{2}^*\wedge f_{3}^*,\gamma\neq
0$&$+\gamma^2(-(f_{1}^*)^2-(f_{2}^*)^2-(f_{3}^*)^2)$&\\

\hline

4.&$f_{1}^*\wedge e_{2}^*\wedge e_{3}^*+ f_{2}^*\wedge
e_{1}^*\wedge e_{3}^*+ f_{3}^*\wedge e_{1}^*\wedge e_{2}^*$&
$(e_{1}^*)^2-(e_{2}^*)^2+(e_{3}^*)^2$&$0$\\

\hline

5.&$f_{1}^*\wedge e_{2}^*\wedge e_{3}^*- f_{2}^*\wedge
e_{1}^*\wedge e_{3}^*+ f_{3}^*\wedge e_{1}^*\wedge e_{2}^*$&
$-(e_{1}^*)^2-(e_{2}^*)^2-(e_{3}^*)^2$&$0$\\

\hline

6.&$f_{3}^*\wedge e_{1}^*\wedge e_{2}^*+ f_{2}^*\wedge
e_{1}^*\wedge e_{3}^*$& $(e_{1}^*)^2$&$0$\\

\hline

7.&$f_{3}^*\wedge e_{1}^*\wedge e_{2}^*- f_{2}^*\wedge
e_{1}^*\wedge e_{3}^*$& $-(e_{1}^*)^2$&$0$\\

\hline

8.&$e_{1}^*\wedge e_{2}^*\wedge e_{3}^*$& $0$&$0$\\

\hline

9.& $0$&$0$&$0$\\

\hline
\end{tabular}
\caption{Classification des $3$-formes effectives en dimension
$6$}{\label{tableau}}
\end{table}

Les invariants $q_\omega$ et $K_\omega$ se d\'eduisent en fait l'un de
l'autre, via la forme symplectique $\Omega$.

\begin{Prop}{\label{compatibilite}}
Soit $\omega\in \Lambda^3_\varepsilon(V^*)$. Pour tout $X,Y\in V$ on a
$$
\Omega(K_\omega X, Y)=q_\omega(X,Y).
$$
\end{Prop}

\begin{proof}
Un calcul direct nous permet de v\'erifier que le r\'esultat est
vrai pour chacune des formes du tableau \ref{tableau}. De plus
pour tout symplectomorphisme $F$ sur $V$ on a
$$
\begin{cases}
q_{F^*\omega}=F^tq_\omega F&\\
K_{F^*\omega}= F^{-1}K_\omega F&\\
\end{cases}
$$
Le r\'esultat est donc vrai pour toute forme.
\end{proof}

\begin{Rem}
Puisque $dq_\omega=q_\omega$, \ref{compatibilite} s'\'enonce
$$
K_\omega=X_{q_\omega},
$$
i.e. $K_\omega$ est le champ hamiltonien associ\'e \`a
l'hamiltonien $q_\omega$. Autrement dit modulo l'identification
des alg\`ebres de Lie $(sp(6),[,])$ et $(S^2(6),\{,\})$ (o\`u
$\{,\}$ est le crochet de Poisson sur $V$ associ\'e \`a $\Omega$),
$q_\omega$ est l'application moment associ\'ee \`a l'action
hamiltonienne de $Sp(6)$. Les diff\'erentes orbites de cette
action sont les diff\'erents niveaux de moment \`a conjugaison
pr\`es.
\end{Rem}

\section{Structures de Monge-Amp\`ere non
d\'eg\'en\'er\'ees et structures g\'eom\'etriques associ\'ees.\\}

\begin{Def}
Une structure de Monge-Amp\` ere $(\Omega,\omega_0)$ sur une vari\'et\'e
$X^6$ de dimension $6$ est dite
\begin{enumerate} 
\item elliptique si $\lambda(\omega_0)<0$,
\item hyperbolique si $\lambda(\omega_0)>0$.
\end{enumerate}
\end{Def}

\begin{Def}{\label{fermee}}
Une structure de Monge-Amp\` ere $(\Omega,\omega_0)$ sur une vari\'et\'e
$X^6$ de dimension $6$ est dite ferm\'ee si la forme normalis\'ee
$$
\omega= \frac{\omega_0}{\sqrt[4]{|\lambda(\omega_0)|}}.
$$
v\'erifie
$$
d\omega=d\hat{\omega}=0
$$
\end{Def}

\subsection{Structures elliptiques.}

\begin{Def}
Une structure presque (pseudo) Calabi-Yau sur une vari\'et\'e r\'e\-el\-le $X^6$
est la donn\'ee 
\begin{enumerate}
\item d'une structure presque (pseudo) K\"ahler $(q,I,\Omega)$ sur $X$:
  $\Omega$ est une $2$-forme non d\'eg\'en\'er\'ee et ferm\'ee, $I$ est une
  structure presque complexe telle $I^*\Omega=\Omega$ et $q$ est la
  m\'etrique (\'eventuellement ind\'efinie) $q=\Omega(I.,.)$;
\item   d'une forme diff\'erentielle complexe $\alpha$ de type $(3,0)$ sur $X$
telle que 
$$
\Omega^3=-\frac{3i}{4}\alpha\wedge\overline{\alpha}.
$$
\end{enumerate}
Une structure presque (pseudo) Calabi-Yau $(q,I,\Omega,\alpha)$ est
une structure (pseudo) Calabi-Yau si la structure presque complexe $I$
est int\'egrable et si la forme $\alpha$ est holomorphe pour cette
structure complexe.
\end{Def}

\begin{Def}
Soit $(\Omega,\omega_0)$ une structure de Monge-Amp\` ere elliptique
sur une vari\'et\'e $X$ de dimension $6$. Soit $\omega$ la forme
normalis\'ee
$$
\omega= \frac{\omega_0}{\sqrt[4]{|\lambda(\omega_0)|}}.
$$
La structure presque (pseudo) Calabi-Yau associ\'ee \`a la classe
conforme $(\Omega,\omega_0)$ est $(q_\omega,K_\omega,\Omega,\alpha)$
avec
$$
\alpha=\frac{4}{3}(\omega+iK_\omega^*(\omega)).
$$
\end{Def}

\begin{Ex}
La structure g\'eom\'etrique sur $T^*\mathbb{R}^3$ associ\'ee \`a
l'\'equation sp\'eciale lagrangienne 
$$
\hess (f)-\Delta f=0
$$
est la structure de Calabi-Yau classique $(q,I,\Omega,\alpha)$ avec:
$$
\begin{cases}
q= -\underset{j=1}{\overset{3}{\sum}} dx_j.dx_j+ dy_j.dy_j&\\
I=\underset{j=1}{\overset{3}{\sum}} \frac{\partial}{\partial y_j}\otimes
dx_j-\frac{\partial}{\partial x_j}\otimes dy_j&\\
\Omega=\underset{j=1}{\overset{3}{\sum}} dx_j\wedge dy_j&\\
\alpha=dz_1\wedge dz_2\wedge dz_3&\\
\end{cases}
$$
De mani\`ere similaire, la structure g\'eom\'etrique associ\'ee \`a
l'\'equation pseudo sp\'eciale lagrangienne 
$$
\hess (f)+\square f=0
$$
est la structure pseudo Calabi-Yau $(q_-,I_-,\Omega,\alpha)$ avec:
$$
\begin{cases}
q_-=dx_1.dx_1-dx_2.dx_2+dx_3.dx_3 +dy_1.dy_1-dy_2.dy_2+dx_3.dx_3&\\
I_-=-\frac{\partial}{\partial y_1}\otimes
dx_1+\frac{\partial}{\partial x_1}\otimes dy_1+ \frac{\partial}{\partial y_2}\otimes
dx_2-\frac{\partial}{\partial x_2}\otimes dy_2-\frac{\partial}{\partial y_3}\otimes
dx_3+\frac{\partial}{\partial x_3}\otimes dy_3&\\
\Omega=\underset{j=1}{\overset{3}{\sum}} dx_j\wedge dy_j&\\
\alpha=dz_1\wedge dz_2\wedge dz_3&\\
\end{cases}
$$
Ces deux structures  sont les seuls exemples associ\'es
\`a des \'equations de Monge-Amp\`ere elliptiques \`a coefficients constants,
puisque, d'apr\`es le tableau \ref{tableau}, il n'y a que deux telles 
\'equations sur $\mathbb{R}^6$.
\end{Ex}

Remarquons qu'\`a toute structure presque Calabi-Yau
$(q,I,\Omega,\alpha)$ est associ\'ee la structure  elliptique
$(\Omega,Re(\alpha)))$. Il y a ainsi une correspondance entre les
structures presques (pseudo) Calabi-Yau et les classes conformes
de structures de Monge-Amp\`ere elliptiques.

\begin{Prop}
La structure presque Calabi-Yau associ\'ee \`a une structure de Monge-Amp\`ere
elliptique est une structure de Calabi-Yau si et seulement si celle-ci est ferm\'ee
au sens de la d\'efinition \ref{fermee}.
\end{Prop}

\begin{proof}
Soit $(\Omega,\omega)$ une structure de Monge-Amp\`ere telle que
$\lambda(\omega)=1$. Soit $I=K_\omega$ la structure presque complexe
associ\'ee et $\alpha=\frac{1}{2}(\omega+iI^*\omega)$ la forme
complexe de type $(n,0)$ associ\'ee. 

Si $(q_\omega,I,\Omega,\alpha)$ est une structure de Calabi-Yau alors
$\alpha$ est holomorphe donc $d\alpha=0$ et de la m\^eme fa\c con,
$d\overline{\alpha}=0$. Donc $d\omega=d\hat{\omega}=0$.

R\'eciproquement supposons que $d\omega=d\hat{\omega}=0$. Alors
$d\alpha=d\hat{\alpha}=0$. Soit sur $TX\otimes \mathbb{C}$
la distribution $C:x\mapsto \{X: i_X\alpha_x=0\}$. Une $1$-forme $\xi$
s'annule sur  $C$ si et seulement $\xi\wedge \alpha=0$ et donc si
$\xi$ s'annule sur $C$ alors $d\xi\wedge \alpha=0$. D'apr\`es le
th\'eor\`eme de Frobenius, la distribution $C$ est int\'egrable et
donc, d'apr\`es le th\'eor\`eme de Newlander-Nirenberg, $I$ est
int\'egrable.
\end{proof}

\subsection{Structures hyperboliques.\\}

Par analogie avec le cas elliptique, le cas hyperbolique est associ\'e
\`a une structure r\'eelle (ou produit) qui est l'analogue des
structures de Calabi-Yau. Ceci nous conduit \`a poser la d\'efinition suivante. 

\begin{Def}
Une structure presque pseudo Calabi-Yau r\'eelle sur une vari\'et\'e
r\'eelle de dimension $6$ est la donn\'ee
\begin{enumerate}
\item d'une structure pseudo K\"ahler hyperbolique $(q,S,\Omega)$:
  $\Omega$ est une $2$-forme non d\'eg\'en\'er\'ee et ferm\'ee, $S$ est
  une structure presque produit ($S^2=Id)$ telle que
  $S^*\Omega=-\Omega$ et $q=\Omega(S.,.)$ est une m\'etrique de
  signature $(3,3)$;
\item de deux formes d\'ecomposables $\alpha$ et $\beta$ dont les
  distributions associ\'ees sont les distributions des vecteurs
  propres de $S$ et telles que
$$
\alpha\wedge \beta= -\frac{\Omega^3}{6}
$$
\end{enumerate}
Une structure presque pseudo Calabi-Yau r\'eelle est une structure pseudo
Calabi-Yau r\'eelle si $S$ est int\'egrable et si $\alpha$ et
$\beta$ sont ferm\'ees.
\end{Def}

\begin{Rem}
Ces vari\'et\'es sont l'analogue des ``vari\'et\'es de
Monge-Amp\`ere'' au sens de Kontsevich et Soibelman (\cite{KS}). Une
vari\'et\'e de Monge-Amp\`ere est pour eux une vari\'et\'e
riemannienne affine $(M,g)$ telle que localement la m\'etrique s'\'ecrit
$$
g= \sum_{i,j} \frac{\partial^2 K}{\partial x_i\partial x_j} dx_i.dx_j,
$$
o\`u $K$ est une fonction lisse v\'erifiant 
$$
det(\frac{\partial^2 K}{\partial x_i\partial x_j})=constant.
$$ 
Comme on le verra plus loin, dans le cas des vari\'et\'es de
Calabi-Yau r\'eelles, on a un tel potentiel $K$ et la m\'etrique s'\'ecrit
$$
g= \sum_{i,j} \frac{\partial^2 K}{\partial x_i\partial y_j} dx_i.dy_j.
$$
avec
$$
det(\frac{\partial^2 K}{\partial x_i\partial y_j})=constant.
$$ 
\end{Rem}

\begin{Def}
Soit $(\Omega,\omega_0)$ une structure de Monge-Amp\` ere hyperbolique
sur une vari\'et\'e $X$ de dimension $6$. Soit $\omega$ la forme
normalis\'ee
$$
\omega= \frac{\omega_0}{\sqrt[4]{\lambda(\omega_0)}}.
$$
La structure presque pseudo Calabi-Yau associ\'ee \`a la classe
conforme $(\Omega,\omega_0)$ est $(q_\omega,K_\omega,\Omega,\alpha,\beta)$
avec
$$
\alpha=\frac{1}{2}(\omega+K_\omega^*(\omega))
$$
et 
$$
\beta=\frac{1}{2}(\omega-K_\omega^*(\omega)).
$$
\end{Def}

\begin{Ex}
La structure (pseudo) Calabi-Yau r\'eelle associ\'ee \`a l'\'equation de
Mon\-ge-Amp\`ere classique
$$
\hess (f)=1
$$
est la structure $(q,S,\Omega,\alpha,\beta)$ avec:
$$
\begin{cases}
q=\underset{j=1}{\overset{3}{\sum}} dx_j.dy_j&\\
S= \underset{j=1}{\overset{3}{\sum}} \frac{\partial}{\partial
  x_j}\otimes dx_j-\frac{\partial}{\partial y_j}\otimes dy_j&\\
\Omega=\underset{j=1}{\overset{3}{\sum}} dx_j\wedge dy_j&\\
\alpha=dx_1\wedge dx_2\wedge dx_3,\;\beta= dy_1\wedge dy_2\wedge dy_3&\\
\end{cases}
$$
\end{Ex}

De fa\c con tout \`a fait similaire au cas elliptique on a la proposition

\begin{Prop}
La structure presque pseudo Calabi-Yau r\'eelle associ\'ee \`a une
structure de Monge-Amp\`ere hyperbolique est une structure pseudo
Calabi-Yau r\'eelle si et seulement si celle-ci est ferm\'ee
au sens de la d\'efinition \ref{fermee}.
\end{Prop}

\subsection{Structures de Monge-Amp\`ere non d\'eg\'en\'er\'ees
localement constantes sur $\mathbb{R}^6$}

\begin{Def}
Une structure de Monge-Amp\`ere $(\Omega,\omega)$ sur une vari\'et\'e
$X$ est dite localement constante si au voisinage de tout point il
existe un syst\`eme de coordonn\'ees de Darboux de $(X,\Omega)$ dans
lequel $\omega$ est \`a coefficients constants
\end{Def}

D'apr\`es $\ref{tableau}$, une structure de Monge-Amp\`ere non d\'eg\'en\'er\'ee localement constante est
associ\'ee \`a une \'equation diff\'erentielle qui modulo un changement de
variable symplectique a l'une des trois formes suivantes:
$$
\begin{cases}
\hess (f)=\gamma,&\gamma\neq 0\\
\Delta f -\gamma \hess (f)=0,&\gamma \neq 0\\
\square f + \gamma \hess (f)=0,&\gamma \neq 0\\
\end{cases}
$$

\begin{Theo}{\label{theo}}
Soit $(\Omega,\omega_0)$ une structure de Monge-Amp\`ere locale
non d\'e\-g\'e\-n\'e\-r\'ee sur une vari\'et\'e r\'eelle de
dimension $6$ $N^6$. Soit la forme normalis\'ee associ\'ee \`a
$\omega_0$
$$
\omega= \frac{\omega_0}{\sqrt[4]{|\lambda(\omega_0)|}}.
$$
$(\Omega,\omega)$ est localement constante si et seulement si
$$
\begin{cases}
d\omega=d\hat{\omega}=0&\\
R(q_\omega)=0&\\
\end{cases}
$$
o\`u $R(q_\omega)$ est le tenseur de courbure de la (pseudo)
m\'etrique $q_\omega$ associ\'e \`a $\omega$.
\end{Theo}

\begin{proof}

\emph{Cas hyperbolique: $\lambda(\omega)=1$\\}

Si $\omega= du_1\wedge du_2\wedge du_3+ dv_1\wedge dv_2\wedge
dv_3$ dans un syst\`eme de coordonn\'ees symplectique $(u,v)$
alors $\hat{\omega}=du_1\wedge du_2\wedge du_3-dv_1\wedge
dv_2\wedge dv_3$. Donc $d\omega=d\hat{\omega}=0$. De plus:
$$
\begin{cases}
K_\omega(\frac{\partial}{\partial u_j})=\frac{\partial}{\partial
u_j}& j=1,2,3\\
K_\omega(\frac{\partial}{\partial v_j})=-\frac{\partial}{\partial
v_j}& j=1,2,3\\
\end{cases}
$$
et donc $q_\omega=\underset{j=1}{\overset{3}{\sum}} du_j.dv_j$:
$q_\omega$ est plate.

R\'eciproquement supposons que $d\omega=d\hat{\omega}=0$ et
$R(q_\omega)=0$. Le probl\`eme \'etant local on peut supposer que
l'on se trouve au voisinage de $a_0=(0,0)$ dans $N=\mathbb{R}^6$.
Soit $\alpha$ et $\beta$ les uniques formes d\'ecomposables telles
que
$$
\begin{cases}
\omega=\alpha+\beta&\\
\hat{\omega}=\alpha-\beta&\\
\alpha\wedge \beta= -\frac{1}{6}\Omega^3&\\
\end{cases}
$$
Puisque $\omega$ et $\hat{\omega}$ sont ferm\'ees, $S=K_\omega$ est
int\'egrable et $\alpha$ et
$\beta$ sont ferm\'ees. Il existe alors un syst\`eme de
coordonn\'ees $(x,y)$ de $N$ dans lequel
$$
\begin{cases}
\omega=f(x,y)dx_1\wedge dx_2\wedge dx_3+ dy_1\wedge dy_2\wedge dy_3&\\
\hat{\omega}=g(x,y)dx_1\wedge dx_2\wedge dx_3-dy_1\wedge dy_2\wedge dy_3&\\
\end{cases}
$$
Mais $\alpha$ et $\beta$ sont ferm\'ees et $\alpha\wedge \beta$ est
constante donc $f$ et $g$ sont des constantes que l'on peut supposer
\'egales \`a $1$.
$K_\omega$ s'\'ecrit alors:
$$
\begin{cases}
K_\omega(\frac{\partial}{\partial x_i})&=\frac{\partial}{\partial
  x_i}\\
K_\omega(\frac{\partial}{\partial y_i})&=-\frac{\partial}{\partial
  y_i}\\
\end{cases}
$$

La forme $\omega$ est effective et en tout point $a$,
$\hat{\omega}_a= K_{\omega_a}^*\omega_a$ donc $\hat{\omega}$ est
effective puisque en tout point $a$
$K_{\omega_a}^*\Omega_a=-\Omega_a$. D\`es lors $\alpha$ et $\beta$
sont effectives. Autrement dit en tout point $a$, $E_\alpha(a)=Ker(\alpha_a)$ et
$E_\beta(a)=Ker(\beta_a)$ sont deux sous espaces lagrangiens de
$(T_aN,\Omega_a)$. En effet si $X$ et $Y\in E_\alpha(a)$ par
exemple:
$$
0= i_{X\wedge Y}(\Omega\wedge \alpha)= \Omega(X,Y)\alpha
$$
D\`es lors:
$$
\Omega= \sum_{i,j=1}^3 a_{ij} dx_i\wedge dy_j
$$

Or $\Omega$ est ferm\'ee. Donc
$$
\begin{cases}
\frac{\partial a_{ij}}{\partial x_k}= \frac{\partial
a_{kj}}{\partial
    x_i}&\\
\frac{\partial a_{ij}}{\partial y_k}= \frac{\partial
a_{ik}}{\partial
    y_j}&\\
\end{cases}
$$
D\'efinissons pour $j=1,2,3$ les fonctions
$$
\phi_j(x,y)=\int_0^{x_1} a_{1j}(t,x_2,x_3,y)dt+ \int_0^{x_2}
a_{2j}(0,t,x_3,y)dt+\int_0^{x_3} a_{3j}(0,0,t,y)dt
$$
Pour $i=1,2,3$ on a $\frac{\partial \phi_j}{\partial x_i}=a_{ij}$.
De plus:
$$
\frac{\partial \phi_j}{\partial y_k}= \frac{\partial
\phi_{k}}{\partial
    y_j}
$$
donc il existe une fonction lisse $\phi$ sur $\mathbb{R}^6$ telle
que $\phi_j=\frac{\partial \phi}{\partial y_j}$ et
$a_{ij}=\frac{\partial^2 \phi}{\partial x_i\partial y_j}$.
Finalement:
$$
\begin{cases}
\Omega=\underset{i,j}{\sum} \frac{\partial^2\phi}{\partial
x_i\partial y_j}
dx_i\wedge dy_j&\\
q_\omega=\underset{i,j}{\sum} \frac{\partial^2\phi}{\partial
x_i\partial y_j}
dx_i . dy_j&\\
\end{cases}
$$

Utilisons maintenant le fait que $q_\omega$ est plate. Soit
$\Gamma_{jk}^j$ les symboles de Christoffels de la connexion de
Levi-Civita $\nabla=\nabla_\omega$ dans les coordonn\'ees $(x,y)$.
Soit pour $j=1\ldots 6$ la matrice
$$
\Gamma_j =(\Gamma_{jk}^l)_{j,k=1}^6
$$
On v\'erifie facilement que pour $j=1,2,3$:
$$
\begin{cases}
\Gamma_j=\begin{pmatrix} \frac{\partial A}{\partial
x_j}A^{-1}&0\\0&0\\\end{pmatrix}&\\
\Gamma_{j+3}=\begin{pmatrix} 0&0\\0&(A^{-1}\frac{\partial
A}{\partial y_j})^{t}\\\end{pmatrix}&\\
\end{cases}
$$
avec $A_{jk} =\frac{\partial^2 \phi}{\partial x_j\partial y_k}$.

 Posons $C_j=\frac{\partial A}{\partial x_j}A^{-1}$ et
$D_j=A^{-1}\frac{\partial A}{\partial y_j}$ pour $j=1,2,3$. De
$R=0$ il vient pour $j,k=1,2,3$:
$$
\begin{cases}
\frac{\partial C_j}{\partial y_k}=0&\\
\frac{\partial C_j}{\partial x_k}-\frac{\partial C_k}{\partial
x_j}+[C_j,C_k]=0&\\
\end{cases}
$$
et
$$
\begin{cases}
\frac{\partial D_j}{\partial x_i}=0&\\
\frac{\partial D_i}{\partial y_j}-\frac{\partial D_j}{\partial
y_i}+ [D_i,D_j]&\\
\end{cases}
$$

Soit $G=G(x)$ une solution du syst\`eme  diff\'erentiel
$$
\frac{\partial G}{\partial x_j} G^{-1}=C_j, \text{ $j=1,2,3$}
$$
(une telle solution existe toujours, c.f. lemme \ref{equa-diff}).
On a alors pour $j=1,2,3$ :
$$
\begin{aligned}
\frac{\partial A^{-1}G}{\partial x_j}&= \frac{\partial
A^{-1}}{\partial x_j}G+ A^{-1}\frac{\partial G}{\partial x_j}\\
& = -A^{-1}\frac{\partial A}{\partial
x_j}A^{-1}G+A^{-1}\frac{\partial G}{\partial x_j}\\
& = -A^{-1}C_jG+A^{-1}C_jG\\
&=0\\
\end{aligned}
$$
et donc:
$$
A(x,y)=G(x)F(y)
$$
Notons $G_j(x)=(G_{j1}(x),G_{j2}(x),G_{j3}(x))$ et
$F_{j}(y)=(F_{1j}(y),F_{2j}(y),F_{3j}(y))$ pour $j=1,2,3$. De
$\frac{\partial^3\phi}{\partial x_j\partial x_k\partial
y_l\phi}=\frac{\partial^3 \phi}{\partial x_k\partial x_j\partial
y_l}$ il vient
$$
<\frac{\partial G_j}{\partial x_k}-\frac{\partial{G_k}}{\partial
x_j},F_l>=0
$$
et donc $\frac{\partial G_j}{\partial x_k}=\frac{\partial
G_k}{\partial x_j}$: il existe des fonctions
$u_1(x),u_2(x),u_3(x)$ telles que
$$
G_j=(\frac{\partial u_{j}}{\partial x_1},\frac{\partial
u_j}{\partial x_2},\frac{\partial u_j}{\partial x_3})
$$
et de la m\^eme fa\c con il existe des fonctions
$v_1(y),v_2(y),v_3(y)$ telles que
$$
F_j=(\frac{\partial v_{j}}{\partial y_1},\frac{\partial
v_j}{\partial y_2},\frac{\partial v_j}{\partial y_3})
$$
D\`es lors dans le syst\`eme de coordonn\'ees $(u,v)$ on a
$$
\begin{cases}
\Omega= du_1\wedge dv_1+ du_2\wedge dv_2+ du_3\wedge dv_3&\\
\omega= r(u)du_1\wedge du_2\wedge du_3+s(v)dv_1\wedge dv_2\wedge
dv_3&\\
\end{cases}
$$
Mais $\omega\wedge \hat{\omega}= -\frac{1}{6}\Omega^3$ donc $r$ et
$s$ sont constantes et inverses l'une de l'autre. Quitte \`a
remplacer $u_1$ par $r u_1$ et $v_1$ par $\frac{1}{r}v_1$ on
obtient le r\'esultat souhait\'e.\\

\emph{Cas elliptique d\'efini n\'egatif: $\lambda=1$, $s(q_\omega)=(0,6)$\\}

La d\'emonstration est analogue \`a la pr\'ec\'edente en rempla\c
cant les coordonn\'ees r\'eelles $(x,y)$ par les coordonn\'ees
complexes $(z,\overline{z})$:

Si
$$
\begin{cases}
\omega=Re(idz_1\wedge dz_2\wedge dz_3)&\\
\Omega= \frac{i}{2}(dz_1\wedge \overline{dz_1}+dz_2\wedge
\overline{dz_2}+dz_3\wedge \overline{dz_3})&\\
\end{cases}
$$
alors $d\omega=d\hat{\omega}=0$ et
$q_\omega=\underset{j=1}{\overset{3}{\sum}} dz_j.d\overline{z_j}$
est plate.

R\'eciproquement supposons que $d\omega=d\hat{\omega}=0$ et
$R(q_\omega)=0$. Comme pr\'ec\'emment le fait que
$\omega=\alpha+\overline{\alpha}$ et $\hat{\omega}$ soient
ferm\'ees implique que la structure complexe $K_\omega$ est
int\'egrable. Il existe un syst\`eme de coordonn\'ees
$(z_1,z_2,z_3)$ de $N$ dans laquelle $\alpha=idz_1\wedge
dz_2\wedge dz_3$. De plus $(q_\omega,K_\omega,\Omega)$ est une
structure de K\"ahler donc il existe un potentiel de K\"ahler
$\phi$ tel que
$$
\Omega=i\partial\overline{\partial} \phi=i\sum_{j,k=1}^3
\frac{\partial ^2\phi}{\partial z_j\partial \overline{z_k}}
dz_j\wedge d\overline{z_k}
$$
De $R(q_\omega)=0$ on d\'eduit que la matrice
$A=(\frac{\partial^2\phi}{\partial z_j\partial \overline{z_k}})$
s'\'ecrit n\'ecessairement
$$
A(z,\overline{z})= G(z)F(\overline{z})
$$
Mais $\overline{A}^t=A$ donc
$\overline{F}^t=GF(\overline{G}^t)^{-1}$. Or $\overline{F}$ et $G$
sont holomorphes et $F(\overline{G}^t)^{-1}$ est antiholomorphe donc constante:
$$
A(z,\overline{z})=H(z)\overline{H}^t(\overline{z})
$$
et donc n\'ecessairement il existe des fonctions holomorphes
$u_1(z),u_2(z),u_3(z)$ telles que
$$
\frac{\partial^2\phi}{\partial z_j\partial
\overline{z_k}}=\sum_{l=1}^3 \frac{\partial u_l}{\partial
z_j}.\frac{\partial \overline{u_l}}{\partial \overline{z_k}}
$$
D\`es lors dans le syst\`eme de coordonn\'ees complexes
$(u_1,u_2,u_3)$ on a le r\'esultat voulu.\\

\emph{Cas elliptique ind\'efini: $\lambda(\omega)=1$,
$s(q_\omega)=(4,2)$\\}

On remarque que la structure de Monge-Amp\`ere pseudo-K\"ahler se
d\'eduit de la structure de Monge-Amp\`ere K\"ahler par le
changement de variable $\xi: z_2\mapsto \overline{z_2}$.
$(\Omega,\omega)$ est alors localement constante si et seulement
si $(\xi^*\Omega,\xi^*\omega)$ est localement constante.

\end{proof}

Nous pouvons r\'esumer la correspondance entre structures (pseudo)
Calabi-Yau et classes conformes de structures de Monge-Amp\`ere
elliptiques dans le tableau \ref{table2}:

\begin{table}[hbp!]
\begin{tabular}{|c|c|}
\hline
presque (pseudo) CY& MA elliptique\\
\hline
(pseudo) CY& MA elliptique ferm\'e\\
\hline
(pseudo) CY plat& MA elliptique localement constant\\
\hline
\end{tabular}
\caption{Correspondance entre structures de Calabi-Yau et structures
  de Monge-Amp\`ere elliptiques}{\label{table2}}
\end{table}

\begin{Ex}
Il existe peu d'exemples explicites de structures de
Calabi-Yau. Le premier exemple non trivial est l'exemple des
m\'etriques de Stenzel sur $T^*S^n$ (\cite{St}) comme  l'explique
Mich\`ele Audin dans son cours \emph{Lagrangian submanifold}
(\cite{A}). C'est un exemple de structure de Calabi-Yau non plate,
aussi l'\'equation sp\'eciale lagrangienne associ\'ee n'est pas
l'\'equation sp\'eciale lagrangienne classique. 

$T^*S^n=\{(u,v)\in \mathbb{R}^{n+1}\times \mathbb{R}^{n+1}: \|u\|=1,
<u,v>=0\}$ s'identifie \`a la quadrique affine complexe $Q^n=\{z\in
C^{n+1}: z_1^2+\ldots + z_{n+1}^2=1\}$ via l'isomorphisme
$$
\xi(x+iy)= (\frac{x}{\sqrt{1+\|y\|^2}},y)
$$
La forme holomorphe sur $Q^n$ peut \^etre d\'efinie par:
$$
\alpha_z(Z_1,\ldots,Z_n)=det_{\mathbb{C}}(z,Z_1,\ldots,Z_n)
$$
Ainsi par exemple dans l'ouvert de carte $z_{n+1}\neq 0$,
$$
\alpha= \frac{(-1)^n}{z_{n+1}} dz_1\wedge \ldots \wedge dz_{n}
$$
La forme de K\"ahler $\Omega= i\partial \bar{\partial} \phi$ est
d\'efinie  \`a partir du potentiel de K\"ahler $\phi=f(\tau)$ o\`u $\tau$ est
la restriction \`a  $Q^n$ de $|z_1|^2+ \ldots + |z_{n}^2$ et $f$ est
solution de l'\'equation diff\'erentielle ordinaire
$$
x(f')^n+ f''(f')^{n-1}(x^2-1)=c>0
$$

Int\'eressons nous au cas $n=3$.  Pour associer \`a une structure de
Monge-Amp\`ere une \'equation ``explicite'' il faut d\'eterminer tout
d'abord un syst\`eme de coordonn\'ees de Darboux. Remarquons pour cela que
$$
\Omega=-dIm\overline{\partial} \phi= d\{f'(\tau)\underset{k=1}{\overset{4}{\sum}}
x_kdy_k-y_kdx_k\}= -2d\{f'(\tau)\underset{k=1}{\overset{4}{\sum}}y_kdx_k\},
$$
puisque sur $Q^3$, $\sum x_kdy_k+ \sum y_kdx_k=0$. D\`es lors dans les
coordonn\'ees $(u,v)$ de $T^*S^3$ on a:
$$
\Omega= -2d\{f'(2+2\|v\|^2)\sqrt{1+\|v\|^2}\sum_{k=1}^4 v_kdu_k\}.
$$
De plus sur l'ouvert de carte $u_4\neq 0$,$du_4=
-\frac{1}{u_4}\underset{k=1}{\overset{3}{\sum}} u_k du_k$ et donc
$\Omega= \underset{k=1}{\overset{3}{\sum}} dw_k\wedge du_k$ avec
$$
w_k =2\frac{f'(2+2\|v\|^2)\sqrt{1+\|v\|^2}}{u_4}(u_kv_4-v_ku_4).
$$
Notons $\psi$ l'application $(u,w)\mapsto (x+iy)$. L'\'equation
sp\'eciale lagrangienne de Stenzel peut s'\'ecrire
$$
(\psi\circ df)^*\omega=0.
$$

\end{Ex}

\section{Annexe}

\begin{Lem}{\label{equa-diff}}
Soit $C_1,C_2,C_3$ des matrices carr\'ees \`a coefficients lisses
au voisinage de $0$ sur $\mathbb{R}^3$ telles que
$$
\frac{\partial C_i}{\partial x_j}-\frac{\partial C_j}{\partial
x_i}+ [C_i,C_j]=0
$$
Alors il existe toujours une solution $G=G(x_1,x_2,x_3)$ solution
du syst\`eme diff\'erentiel
$$
\begin{cases}
\frac{\partial G}{\partial x_1}G^{-1}=C_1&\\
\frac{\partial G}{\partial x_2}G^{-1}=C_2&\\
\frac{\partial G}{\partial x_3}G^{-1}=C_3&\\
\end{cases}
$$
\end{Lem}

\begin{proof}
Soit $X(x_1,x_2,x_3)$ l'unique solution de
\begin{equation}{\label{1}}
\frac{\partial G}{\partial x_1}G^{-1}=C_1
\end{equation}
avec $X(0,x_2,x_3)=Id$. Pour n'importe quelle matrice
$Y(x_2,x_3)$, $XY$ est encore solution de $\eqref{1}$. Nous allons
montrer que l'on peut choisir $Y(x_2,x_3)$ telle que $XY$ soit
solution du syst\`eme.

$XY$ est solution du syst\`eme si et seulement si
\begin{equation}{\label{sys2}}
\begin{cases}
\frac{\partial Y}{\partial x_2} Y^{-1}=X^{-1}(C_2X-\frac{\partial
X}{\partial x_2})&\\
\frac{\partial Y}{\partial x_3} Y^{-1}=X^{-1}(C_3X-\frac{\partial
X}{\partial x_3})&\\
\end{cases}
\end{equation}

Posons $C'_2=X^{-1}(C_2X-\frac{\partial X}{\partial x_2})$ et
$C'_3=X^{-1}(C_3X-\frac{\partial X}{\partial x_3})$. On v\'erifie
sans peine que
$$
\begin{cases}
\frac{\partial C'_2}{\partial x_1}=0&\\
\frac{\partial C'_3}{\partial x_1}=0&\\
\frac{\partial C'_2}{\partial x_3}-\frac{\partial C'_3}{\partial
x_2}+ [C'_2,C'_3]=0&\\
\end{cases}
$$
Soit alors $\tilde{Y}(x_2,x_3)$ l'unique solution de
\begin{equation}{\label{2}}
\frac{\partial \tilde{Y}}{\partial x_2}\tilde{Y}^{-1}=C'_2
\end{equation}
telle que $\tilde{Y}(0,x_3)=Id$. $\tilde{Y}(x_2,x_3)Z(x_3)$ est
solution de $\eqref{sys2}$ si et seulement si $Z(x_3)$ est
solution de
\begin{equation}{\label{3}}
\frac{\partial Z}{\partial x_3}Z^{-1}=C"_3
\end{equation}
avec $C"_3=\tilde{Y}^{-1}(C'_3\tilde{Y}-\frac{\partial
\tilde{Y}}{\partial x_3})$. Or on v\'erifie que
$$
\frac{\partial C"_3}{\partial x_2}=0
$$
Donc il existe un unique $Z(x_3)$ solution de $\eqref{3}$ tel que
$Z(0)=Id$.  Une solution du syst\`eme initial est alors:
$$
G(x_1,x_2,x_3)=X(x_1,x_2,x_3)\tilde{Y}(x_2,x_3)Z(x_3)
$$
\end{proof}

\end{document}